\documentclass [11pt,a4paper]{article}

\usepackage{amsmath} 
\usepackage{amsthm}
\usepackage{graphics}
\usepackage{latexsym} 
\usepackage{amssymb} 
\usepackage{enumerate}

\newtheorem{lemma}{Lemma}

\newtheorem{proposition}{Proposition}

\newtheorem*{main}{Theorem 1}
\def \N {\mathbb N}
\def \Z {\mathbb Z}

\def \R {\mathbb R}

\def \ind{1\!\!1}
\def \vareps {\varepsilon}

\def \expect      {\ensuremath{\mathbf{E}}}

\begin{document}

\title{A signal-recovery system: asymptotic properties, and construction
of an infinite-volume limit}

\author{J. van den Berg and B. T\'{o}th \\
CWI Amsterdam and TU Budapest
}
\date{}

\maketitle

\begin{abstract}
We consider a linear sequence of `nodes', each of which can be in state
$0$ (`off') or $1$ (`on'). Signals  from outside are sent to the
rightmost node and travel instantaneously as far as possible to
the left along nodes which are `on'. These nodes are immediately
switched off, and become on again after a recovery time. The
recovery times are independent exponentially distributed  random
variables.  

We present properties for finite systems and use some of these properties to 
construct an infinite-volume extension, with signals `coming from infinity'.
This construction is related to a question by D. Aldous and we expect that
it sheds some light on, and stimulates further investigation of, that question.
\end{abstract}

\newpage

\section{Introduction}

Let $X_1(t), \ldots X_n(t)$ be $0$--$1$-valued random processes
described as follows: When $X_i$ equals $0$ it becomes 1 at rate $\rho_i$,
independently of the other $X_j$-s. If each of $X_i, X_{i+1}, \ldots,
X_n$ equals 1, then at rate $\rho$ they all become simultaneously $0$. 
We start at time $t=0$ with all $X_i's$ equal to $0$.
The $\rho_i$-s and $\rho$ are parameters of the model, called
\emph{recovery rates} and the \emph{input rate}, respectively.

This system can be interpreted as a simple model of a communication
line, and we will frequently use terminology motivated by this
interpretation: the indices $1,2, \ldots, n$ correspond to 
nodes which can be `on' (have value 1) or `off' (have value 0).
Signals from outside are sent at rate $\rho$ to the rightmost 
node $n$ and are  transmitted
instantaneously as far as possible to the left until they are blocked by 
an off-node. The nodes passed by the signal are switched off
immediately. When a node $i$ is `off', it becomes `on' 
after an exponentially distributed (with mean $1/\rho_i$)
recovery time. Recovery times are completely independent. 

Another interpretation is in terms of forest fires (or infections):
the numbers $1, \ldots, n$ are possible locations of trees. At the
rightmost location
ignition attempts are made at rate $\rho$, and an attempt succeeds if
that position is occupied. When a tree is on fire, it
immediately sets the tree on the next location to
its left (if that location is
occupied)
on fire and disappears (i.e. its position becomes empty). If position $i$
is empty, a new tree appears there at rate $\rho_i$.
Since in a real forest
the growth of new trees is much slower than the propagation of fire, the
instantaneous spread of fire (or infection) in our model is not as unrealistic
as it may seem at first sight. (The one-dimensionality {\em is} of course
a serious simplification in this context).

In the above description the incoming signals correspond to a Poisson process.
More generally we will consider signals corresponding to a
renewal process. The distribution function of the intervals between
consecutive incoming signals will be denoted by $F^{(n+1)}$.  
(The reason for this notation, with the superscript $n+1$, will
become clear later).

So, more precisely, this more general model is as follows: The parameters
of the model are $F^{(n+1)}$ (a distribution function with $F^{(n+1)}(0) = 0$),
and the recovery rates
$\rho_i$, $i =1, \ldots, n$. Introduce i.i.d. random variables 
$\tau_1, \tau_2, \ldots$ with distribution function $F^{(n+1)}$, and call 
the values $T_i:=\sum_{j=1}^i \tau_j$, $i=0,1, 2, \ldots$ 
\emph{input times}.
At the zeroth input time $T_0=0$ we set each $X_i$, $i=1, \ldots, n$
equal to $0$. When $X_i = 0$ it becomes $1$ at rate $\rho_i$,
independent of the other $X_j$-s and of the $\tau_j$-s.
If, at time $t$, $X_i = X_{i+1} = \ldots, =X_n = 1$, 
then each $X_j$, $j \geq i$
becomes 0 at the smallest input time larger than $t$.
We call this model a 
\emph{(size $n$) on-off system with recovery rates
$\rho_1, \ldots, \rho_n$ and input interval distribution 
function $F^{(n+1)}$}. The case mentioned in the first paragraph, 
when the input signals arrive according to a Poisson process of
intensity $\rho$, corresponds to $F^{(n+1)}=E^{\rho}$, where
$E^\rho(t)=1-e^{-\rho t}$ denotes the exponential distribution
function with expectation $1/\rho$. 

As said before, we will frequently use terminology inspired by the 
signal interpretation. Although this terminology is practically
self-explanatory, we want to define some of these terms more precisely,
to avoid confusion: we say that a signal \emph{is sent to node
$i$} at time $t$, 
if $X_{i+1}$ switches from on to off at time $t$ (or, in case $i$ is the
rightmost node, if $t$ is an input time). And we say that a signal
\emph{is received by $i$} at time $t$, if $X_{i}$ itself switches from 
on to off at time $t$.

Finally, we will also consider the case that input signals are generated
`permanently'.
By this we mean that the rightmost node, $n$, after each recovery
immediately 
receives a signal (and hence is switched off again). In this
case 
we say (with some abuse of notation, since there are no proper input
intervals anymore) that $F^{(n+1)}=[0]$. 
It is easy to see that
this case is (when we only `observe' the behaviour of the
nodes $1, \ldots, n-1$) equivalent to the earlier mentioned case
with $n-1$ nodes and with
Poisson (intensity $\rho_n$) input signals, i.e., with input
interval distribution function $F^{(n)}=E^{\rho_n}$.

Several interesting questions arise: Suppose the input is
Poissonian, and all 
recovery rates are equal
(say 1), and we start with all nodes empty. What is the
asymptotic behaviour 
(as $n \to  \infty$) of the expectation of the first time a
signal arrives at node $0$. This appears to be of order $\log
n$. 
(Note that this means that when the system is in equilibrium, the
probability 
that an arbitrary incoming signal reaches node $0$, is of order
$1/(\log n)$, 
which is much larger than one would naively guess). 
This is done in Section 2. Several
arguments in that section are of key importance for Section 3,
which deals with 
the question whether there are non-trivial extensions of this
model 
to \emph{infinite systems, with signals ``coming from
infinity"}. 
The answer, as stated in Theorem 1,  is positive and is
related to a question posed by 
D. Aldous.  We hope the result
sheds some more light on that question. The proof of Theorem 1
is presented in Section 4. 

\section{Properties of the finite system}

Consider a size $n$ on-off system (as defined in Section 1)
with input interval distribution function $F^{(n+1)}$ and
recovery rates  $\rho_1, \ldots, \rho_n$.
As stated before, the input signals correspond to a renewal
process. 
It is easy to see that the times at
which signals are {\em received} by $n$ (i.e. the times at which
$X_n$ 
switches from ``on" to ``off") also form a renewal process
(because, whenever 
$X_n$ switches from ``on" to ``off", the process, as far as node
$n$ is 
concerned, starts afresh). Since each signal received by $n$ is
sent 
instantaneously to $n-1$, we can repeat the above argument and
conclude 
that the reception times of signals at $n-1$ also form a renewal
process etc. 
We call the distribution function of the difference between
consecutive 
times at which node $i$ receives a signal, its {\em
interreception time  
distribution function}.

The following lemma relates the interreception time  distributions of
two consecutive nodes:

\begin{lemma}
\label{firstlemma}
Let, for $1 \leq i \leq n$,  $F^{(i)}$ be the interreception time
distribution function of node $i$ and $F^{(n+1)}$ the input interval
distribution. 
Define:
$$
\phi^{(i)}(s) = 1 - \int_0^{\infty} e^{-s x} \, dF^{(i)}(x), \,
\quad
i=1,2, \ldots,n+1
$$
Then:
$$
\phi^{(i)}(s) = \frac{\phi^{(i+1)}(s)}{\phi^{(i+1)}(s + \rho_{i})}, 
\quad
i = 1, \ldots, n.
$$
\end{lemma}

\begin{proof}
Let $\tau$ be the first time node $i$ switches from `off' to
`on',  and let $Y$ be the first time it receives a signal. 
Further, let $Z_k$ be the $k$-th time node $i+1$ receives a
signal, and 
let $\xi_{k} = Z_{k}-Z_{k-1}$, $k=1,2,\dots$. 
The random variable $\tau$ is exponentially distributed with
parameter $\rho_i$. Furthermore, the random variables 
$\xi_k$, $k=1,2,\dots$  are i.i.d.   and also 
independent of $\tau$. So  we have:
\begin{align}
\phi^{(i)}(s)=
1-\expect\big(e^{-sY}\big)
&
=
1-
\sum_{k=0}^\infty
\expect\big(
e^{-sZ_{k+1}}\ind_{ \{ \tau \in [Z_{k},Z_{k+1}) \} } 
\big)
\notag
\\
&
=
1-
\sum_{k=0}^\infty
\expect\big(
e^{-sZ_{k+1}}
\big(e^{-\rho_i Z_k} -e^{-\rho_i Z_{k+1}}\big) 
\big)
\notag
\\
&
=
1- 
\sum_{k=0}^\infty
\expect\big(
e^{-(s+\rho_i)Z_{k}}
\big(e^{-s \xi_{k+1}} -e^{-(s+\rho_i) \xi_{k+1}}\big) 
\big)
\notag
\\
&
=
1-
\frac
{\expect\big(e^{-s\xi_1}\big)-
\expect\big(e^{-(s+\rho_i)\xi_1}\big)}
{1-\expect\big(e^{-(s+\rho_i)\xi_1}\big)}
\notag
\\
&
=
\frac{\phi^{(i+1)}(s)}{\phi^{(i+1)}(\rho_i+s)}.
\end{align}

\end{proof}

By repeated application of the above Lemma 1, and using
induction, we get 

\begin{lemma}
\label{secondlemma}
For $1 \leq i \le j \leq n+1$:
\begin{equation} 
\label{mform}
\phi^{(i)}(s) = 
\frac{\prod_{A \subset \{i,\ldots, j-1\}\, : \,  
|A| {\mathrm{\ even}} } \,
\phi^{(j)}(s+ \sum_{k \in A} \rho_k)}
{\prod_{A \subset \{i,\ldots, j-1\} \, : \,  
|A| {\mathrm{\ odd}} } \,
\phi^{(j)}(s+ \sum_{k \in A} \rho_k)},
\end{equation}
where $|A|$ denotes the number of elements of $A$.
\end{lemma}

This immediately gives the following result: 

\begin{lemma}
\label{thirdlemma}
The interreception time distribution of node $i$,  $F^{(i)}$, 
is invariant under permutations of the
sequence of recovery rates $\rho_i, \rho_{i+1}, \ldots, \rho_n$.
\end{lemma}

\noindent
{\sl Remark:} 
In spite of its apparent simplicity, this observation is rather
surprising: it easily follows from identity (\ref{mform}) but we
do not see any simple direct `pathwise' argument  for its proof. 

Lemma 3 is important in the construction of an infinite-volume
system in the next section. We will illustrate its strength in
the remainder   of the present section. 
We consider  the special case when all $\rho_i$-s are equal, say
$1$, and  the inputs come permanently (that is,
$F^{(n+1)}= [0]$). As already mentioned, this is equivalent to a system
of $n-1$ nodes with recovery rates 1 and Poissonian input with
rate 1, so that $\phi^{(n)}(s) = s/(1+s)$.  Using the identity
$\binom{n-1}{l-1} + \binom{n-1}{l} = \binom{n}{l}$, from
(\ref{mform}) we get 
$$
\phi^{(1)}(s) = 
\frac
{\prod_{0 \leq k \leq n \, : \,  k \, 
{\mathrm{\ even}} }
(s+k)^{\binom{n}{k}}}
{\prod_{0 \leq k \leq n \, : \, k \, 
{\mathrm{\ odd}} }
(s+k)^{\binom{n}{k}}}.
$$
We denote here 
by $T_n$ the first time a signal is received by node 1.
As we are interested in the asymptotics for long chains of nodes,
we denote explicitly by the subscript $n$ the length of
the string of identical nodes considered. 
Thus 
$\phi^{(1)}(s)= 1- \expect\big(\exp(-s T_{n})\big)$.
By evaluating the derivative of the above expression at $s=0$, 
we get
\begin{equation}
\label{exp}
\expect\big(T_{n}\big) = 
\frac
{\prod_{1 \leq k \leq n \, : \, k \, 
{\mathrm{\ even}} }
k^{\binom{n}{k}}}
{\prod_{1 \leq k \leq n \, : \, k \, 
{\mathrm{\ odd}} }
k^{\binom{n}{k}}}.
\end{equation}
%
%and, taking logarithm, 
%\begin{equation}
%\label{logexp}
%\log \expect\big(T_{n}\big) = 
%\sum_{k=1}^{n+1} (-1)^k \binom{n+1}{k} \log k.
%\end{equation}

Andr\'as Luk\'acs, \cite{lukacs}, drew our
attention to the survey-article \cite{flajoletsedgewick} 
by Flajolet and Sedgewick, about the use of contour integrals
(and Melling transforms) to study the asymptotic behaviour as
$n \to \infty$ of expressions of the form
$\sum_{k=1}^n (-1)^k \binom{n}{k} f(k)$ for a wide range of functions $f$.
The case $f(k) = \log k$ is one of the examples they handle
(see their Theorem 4), and according
to their paper the expression in the r.h.s. of  (\ref{exp}) is asymptotic
to $e^{\gamma}$. So
\begin{equation}
\label{limit}
\lim_{n\to\infty}
\frac{\expect\big(T_{n}\big)}{\log n}
= e^{\gamma},
\end{equation}
where $\gamma= 0.577\dots$ is Euler's constant.
Although the following quite elementary probabilistic argument,
based on
Lemma 3 above, does not give the precise value of the limit in
(\ref{limit}), 
it does give the correct order of magnitude of
${\expect\big(T_{n}\big)}$ 
One of our reasons for working this out here is that a similar argument
is used in the construction of the infinite-volume system in
Section 4. Another reason is that from the paper by Flajolet and
Sedgewick,  \cite{flajoletsedgewick}, 
one gets the impression that no elementary way is known to obtain the
order of magnitude of ${\expect\big(T_{n}\big)}$.

\begin{proposition}
Consider, for each $n$,  a finite on-off system with nodes
$\{1, \ldots, n\}$, where all
recovery rates are $1$, and with permanent input  signals.
Let $T_n$ denote the first time node 1 receives a signal. Then  
there exist constants $C_1, C_2 > 0$ such that for all $n$
\begin{equation}
\label{bounds}
C_1 < \frac{\expect\big(T_{n}\big)}{\log n} < C_2
\end{equation}
\end{proposition} 

%\noindent
%{\sl Remark:} 
%A refinement of the
%arguments below shows that ...

\begin{proof}
We use stochastic domination in proving both bounds.

The \emph{lower bound} is easy: 
Note that before the first receival time at node 1 all 
nodes $1,2,\dots n$ must recover at least once. So $T_n$
stochastically dominates $\max \{\tau_i : 1 \leq i \leq n\}$, 
where $\tau_1,\tau_2,\dots,\tau_n$ are i.i.d. exponentially
distributed random variables with mean 1. It follows that 
$$
\expect\big(T_n\big)
\ge
\expect\big(\max \{\tau_i : 1 \leq i \leq n\}\big)
=
\sum_{i=1}^n \frac1i
=\log n + {\cal O}(1),
$$
which proves the lower bound. 

The \emph{upper bound} uses a little
trick: Suppose we add
an extra node 0 at the left of node 1, with recovery rate $1 / \log n$.
Denote this new system by II and the old system by I.
Let $\tilde{T}$ be the first time
in system II that node 0 receives a signal. It is clear that system II
is an extension of the old one, in the sense that the nodes
$1, \ldots, n$
`do not feel the change', so that obviously $T \leq \tilde{T}$.
Finally consider the system, denoted III,  obtained from system I by putting
an extra node $n+1$ at the right of $n$, with recovery rate $\rho_{n+1}
= 1/(\log n)$. (So, in system III the input signals are sent to
$n+1$ which, if it is `on',  sends them to $n$, etc).
Let $\hat{T}$  denote the first time node 1 receives a signal in
system III. By Lemma 2, $\hat{T}$ has the same distribution as $\tilde{T}$.
So we have
$$
\expect(T) \leq \expect(\tilde{T}) = \expect(\hat{T}).
$$
The following computation is for system III.
Let $k$ be a non-negative integer.
Let $A$ be the event that an input signal is sent in the time interval
$(4 k \log n, (4 k+ 1) \log n)$, $B$ the event that node $n + 1$ has
no recovery in the interval $(4 k \log n, (4 k+ 2) \log n)$, but does
have a recovery in $((4 k + 2) \log n, (4 k + 3) \log n))$,
$C$ the event that each of the nodes $1, \ldots,n$ which is off at
time $(4 k +1) \log n$ has a recovery before time $(4 k + 2) \log n$,
and $D$ the event that an input signal is sent to $n+1$ in the
interval $((4 k + 3) \log n, (4 k + 4) \log n)$.
It is easy to see that the conditional probability of $A \cap B \cap C \cap D$
given all information up to time $4 k \log n$ is at least
$$ 
(1- e^{- \log n}) \,  e^{-2} \, (1- e^{-1}) \,  (1 - e^{- \log n})^n \,
(1 - e^{- \log n}),
$$
which is larger than  $ \alpha:= e^{-3}(1-e^{-1})/2 > 0$, 
uniformly in $k$, for sufficiently large $n$.
Moreover, if all the events $A$-$D$ happen, node $1$ will receive a signal
in the interval $((4 k+ 3) \log n, (4 k+4) \log n)$ (and hence
in $(4 k \log n, 4 (k+1) \log n)$).
So, for each integer $k \geq 1$, we have $P(\hat{T} > 4 k \, \log n)
\leq \alpha^k$, from which the required result follows.

\end{proof}

%\noindent
%{\sl Remark:} 
%By inserting \emph{two} extra sites in the system, with rates
%$(x\log n)^{-1}$ and $(y\log n)^{-1}$, using the same type of
%argument we get 
%$$
%\limsup_{n\to\infty}\frac{\expect\big(T_n\big)}{\log n}\le 
%\inf_{x,y>0}
%\frac{x^2-y^2}{x\exp(-1/x)-y\exp(-1/y)}=
%2.38105\dots.
%$$
%Inserting more and more extra `filtering' nodes one could get better
%and better upper bounds. 

\section{Infinite-volume models}

Note that a finite on-off system, as introduced in Section 1, could be
described as a collection $X_1(t), \ldots, X_n(t)$ of 
$0$--$1$ valued processes with the property that the time intervals during
which a process has value
 $0$ are
independent, exponentially distributed (those for $X_i$ with mean
$1/ \rho_i$), and that, after independent time intervals with distribution
$F^{(n+1)}$, the string of $1$-s connected to node $n$ is turned into $0$-s.
In this section we investigate the question whether there are
suitable infinite-volume systems with such properties.
There are several cases to distinguish, depending on
the asymptotic behaviour of the $\rho_i$-s and the nature of the 
input signal `at infinity' (which will be made precise later). 
The most interesting appears to be the case where 
\begin{equation}
\label{casefour}
\sum_{i=1}^{\infty}
e^{- t \rho_i} < \infty,
\qquad
\forall t > 0,
\end{equation}
and with `permanent input signals at infinity'.
In the present paper we consider only this case in
detail. However, see Remark (iii) after Theorem 1 stated
below for a concise description of the other possibilities. 

The above condition (\ref{casefour}) on the $\rho_i$-s  means, by
Borel-Cantelli, 
that if we start with all nodes in state $0$, and there would be
NO input signals,
there is an infinite connected string of $1$-s at any positive  
time $t > 0$.
So, when we do take into account permanent input signals at infinity
we expect, intuitively, that in every time interval, no matter
how small,   infinite connected strings of
1's are formed and immediately destroyed (i.e. turned into $0$-s).
It is not at all clear at this
stage that a dynamics with such kind of behaviour exists; see
Remark (i)  below about existence problems for so-called
frozen-percolation models, and  
Remark (ii).

The main result of the present paper is a proof that such a
system does indeed exist. 
More precisely, we prove the following theorem.

\begin{main}
Let $\rho_i, i=1, 2, \ldots$ be positive numbers satisfying
(\ref{casefour})
There exist $0$--$1$-valued processes
%$\R_+\ni t\mapsto X_i(t)\in\{0,1\}$, $i\in\N$,
$ X_i :\R_+ \mapsto \{0,1\}$, $i\in\N$,
defined jointly 
on the same probability space, with the following properties:
\begin{enumerate}[(a)]
\item
Almost surely, for all $i\in\N$,  $X_i(0) = 0$.
\item
Almost surely, for all $i\in\N$,  $t\mapsto X_i(t)$  
is continuous from the right having left limits (c.a.d.l.a.g.).
\item
Let $T_k^i$ denote the length of the $k$th interval during
which $X_i(\cdot)$ equals 0. Then each $T_k^i$ is exponentially
distributed with mean $1/\rho_i$, and the random variables
$\big(T_k^i\big)_{i,k\in\N}$ are independent. 
\item
Almost surely, for all $t\in\R_+$ and $k\in\N$ with 
$X_k(t^-) = 1$:  if for all $l\ge k$ \ 
$X_l(t^-) = 1$ then 
$X_k(t) = 0$, else 
$X_k(t) = 1$.
\end{enumerate}
Moreover, the collection of processes
$t\mapsto X_i(t)$ , $i=1, 2, \ldots$ has the following
additional properties: 
\begin{enumerate}[(a)]
\setcounter{enumi}{4}
\item
Almost surely, there are no $t$ and $k$ such that 
$X_l(t) =1$ for all  $l\ge k$.
\item
Almost surely,  the reception times of signals are dense. That
is, for all $t \geq 0$ and $\varepsilon > 0$ 
there exist $i\in\N$ and $s\in(t,t+\vareps)$,
such that $X_i(s^-) = 1$ and $X_i(s) = 0$.
\end{enumerate}
\end{main}

\noindent
{\bf Remarks:} 

\smallskip\noindent
(i)
The following remark illustrates why the existence of such a
process is not 
obvious:  D. Aldous , \cite{aldous}, has introduced a
percolation model in which infinite clusters are `frozen' (we will
refer to this model as `frozen-percolation').
Informally, 
that model is as follows. Each vertex (or, for
\emph{bond percolation}, each edge) of a countably infinite,
locally finite connected  graph $G$ can have state $0$ or $1$. At
time $0$ they are all in state $0$. 
Now, assign to each vertex $i$
a time $\tau_i$. The $(\tau_i)$ are iid random variables
with a continuous distribution. Each vertex $i$ remains $0$ until
time $\tau_i$. Then it switches to $1$ (and stays 1 forever),
unless some neighbour of $i$ 
already belongs to an infinite cluster of $1$-s, in which case
$i$ remains 0 
forever. Aldous constructed such a process for the case where $G$
is the 
regular binary tree, and posed the question whether it exists for
$\Z^d$. Benjamini and Schramm, \cite{benjaminischramm},
have pointed out that it does not exist  for $\Z^2$. The
following simple, deterministic, 1-dimensional example, 
due to Antal J\'arai, \cite{jarai}, shows very clearly the
essence of the difficulty: 

\begin{proposition} [A. J\'arai, \cite{jarai}]
Let $t_1, t_2, \ldots$ be a sequence of distinct,\\
 strict\-ly
positive numbers which tends to $0$.
There is \emph{no} sequence of functions
$\omega_i : \R_+ \mapsto \{0,1\}$, $i=1,2,\dots$  with
the following properties: 
$$
\omega_i(t) := 
\left\{ 
\begin{array}{ll}
0       & \mbox{ \ if \ }            t < t_i 
          \mbox{ \ or \ }            \omega_j(t_i^-) = 1
          \mbox{ \ for all \ }       j > i
\\
1       & \mbox{ \ otherwise}
\end{array}
\right. 
$$
\end{proposition}

\begin{proof}

Suppose such a sequence does exist. There are two possibilities: either
there exist $t$ and $i$ with $\omega_j(t) = 1$ for all $j \geq i$
or there exist no such $t$ and $i$. In the latter case we have (by
the rules above) that $\omega_j(t) = 1$ for all $j$ and all $t \geq t_j$.
Since all $t_j$ are smaller than some number $t_{\mbox{max}}$, every
$\omega_j$ equals 1 at time $t_{\mbox{max}}$, a contradiction.
As to the former case, let $t$ and $i$ be as stated there. Let $j$
be the smallest number larger than $i$ with $t_k < t_i$ for all
$k \geq j$. From the rules given above (and the assumption for this case)
it follows that
$\omega_k(t_{j-1}^-) = 1$ for all $k \geq j$ and so $\omega_{j-1} =
0$ at {\em every} time, in particular at time $t$: again a contradiction.
Since both cases lead to a contradiction, the proposition has been proved.

\end{proof}

Note that, when the $t_i$ are not deterministic but independent,
exponentially distributed  random variables with
mean $1/\rho_i$, $i=1, 2, \ldots$, with the $(\rho_i)$ as in 
Theorem 1 (i.e., in our terminology, when they are the first recovery times
of the nodes in the system Theorem 1 deals with) the condition in
J\'arai's example is satisfied with probability 1.
This shows that the frozen percolation model does not exist for the half-line
with $\rho_i$-s as in Theorem 1. \\
Although the situation for $\Z^2$ looks more complicated than for the
half-line, the reason why frozen percolation does not exist is essentially
the same:
Consider frozen percolation on  $\Z^2$, with identically (say:
exponentially) distributed holding times $\tau_i$, $i\in\Z^2$.
At the critical time (just before an infinite
cluster is formed), there are infinitely many separate (not connected
with each other) open circuits  around the origin. If we then consider the
sequence $(t_i), i =1,2, \ldots$ of (random) times needed to connect
consecutive circuits, we are exactly in the situation of J\'arai's example.
This illustrates how study of the half-line can give more insight on
what happens on $\Z^2$. 

\smallskip\noindent
{(ii}) 
In Section 5 of his paper, D.  Aldous poses some open questions
related to the 
frozen-percolation model. One of them is whether a dynamics
exists  where vertices (or, for bond percolation, edges) become
open (in state 1)  at rate 1 and where infinite clusters of $1$-s
are destroyed   (i.e. turned into 0's) immediately. Although this
question was formulated  for graphs which have critical
percolation probability less than $1$,  like the binary tree
or $\Z^d$, $d\geq 2$, we think, based on the previous Remark,
that results for analogous problems on the  half-line, like our Theorem 1,
can help to better understand these problems. 

\medskip\noindent
{(iii)}
Now, returning to the setup of the present paper: given the
recovery rates $\rho_i$, $i=1,2,\dots$ define
$$
\theta :=
\sup\{t\in\R_+: \sum_{i=1}^\infty e^{-\rho_it}=\infty\}=
\inf\{t\in\R_+: \sum_{i=1}^\infty e^{-\rho_it}<\infty\}.
$$
There are four essentially different cases with essentially
different behaviour of the infinitely extended system.
Theorem 1 refers to Case 4, the only really interesting
one. The claims below for Cases 1,2 and 3, which are
formulated in a quite informal way, can be stated more precisely,
and proved by
straightforward applications of Borel-Cantelli lemmas.  
\\
{\bf Case 1:} 
If $\theta=\infty$, then by a simple Borel-Cantelli
argument one can see, that in the infinitely extended system
\emph{no signals coming from infinity will penetrate the
system}. This is the case when $\rho_k\ll\log k$, as
$k\to\infty$. The system with constant recovery rates,
$\rho_i=1$, belongs to this case.
\\
{\bf Case 2:} 
If $\theta<\infty$ and 
$\sum_{i=1}^\infty e^{-\rho_i\theta}=\infty$ then one can
construct an infinite dynamics which satisfies
properties (a-c)
stated in Theorem 1, but not properties (d-f) (inclusion of
property d in this case leads to the same kind of problems as in 
Jarai's example (see Remark (i) above)). 
In particular, there will be nonempty time intervals during which infinite
connected strings of $1$s are present in the system. This makes
the dynamics uninteresting. Typical example is
$\rho_k=\theta^{-1}\log k$. 
\\ 
{\bf Case 3:}
If $\theta<\infty$ and 
$\sum_{i=1}^\infty e^{-\rho_i\theta}<\infty$ then one can
construct an infinitely extended dynamics with moderately
interesting behaviour. Namely: in this case, if at some time
$t_0$ all but finitely many nodes are in state $0$, then exactly
at time $t_0+\theta$ an infinite connected string of $1$-s 
emerges (Borel-Cantelli), which is instantaneously erased by a
signal  penetrating from infinity and sweeping
through the system, down to the rightmost node in state $0$.
So, one can construct with `bare hands' a
dynamics where periodically, with period $\theta$, signals
penetrate from infinity and erase an infinite  connected string
of $1$s, 
just emerging. Typical example is 
$\rho_k=\theta^{-1}\log k+\alpha\log\log k$, with $\alpha>1$. 
\\
{\bf Case 4:}
The only really interesting case is $\theta=0$. In this case 
infinite connected strings of $1$-s try to emerge `in no time'
and are immediately swept away by signals penetrating from infinity. So the
constructed system is in a permanent state of excitation. This
behaviour is intuitively somewhat related to the so-called
\emph{self-organized criticality} phenomenon which recives enormous
attention in the
physics literature. This case is the subject of Theorem 1.

\medskip\noindent
(iv)
A very natural question to ask is whether properties (a)-(d)
listed in Theorem 1 determine \emph{uniquely} the process.
Under the extra condition that the signal reception times at each
node form a \emph{renewal process}, we can prove uniqueness. This 
uses very similar ideas to the ones presented in the forthcoming  
proof of Theorem 1. We do not include  this proof in the
present note. We can not anwer this question in full generality, 
without the extra assumption mentioned above.

\section{Proof of Theorem 1}

To prove the main theorem we will first revisit the finite case
studied in Sections 1 and 2, and introduce some more terminology
and notation.  So consider a finite on-off system with nodes
$\{1, \cdots, n\}$, recovery 
rates $\rho_1, \cdots, \rho_n$, and input 
interval distribution function $F^{(n+1)}$. Suppose at time $0$ all nodes
are off. 
Let, for  $1 \leq i \leq n$ and  $k = 1, 2,  \ldots$, $R_k^i$
denote 
the $k$th recovery time at node $i$, i.e., the $k$th time it
switches from 
`off' to `on'. Also,  let $S_k^i$ be the $k$th time a signal is
received by 
node $i$. For convenience, we will define $S_0^i = 0$. 
Let $\mu(F^{(n+1)}; \rho_n, \rho_{n-1}, \cdots, \rho_1)$ denote the joint 
distribution
of the collection 
$\big((R_k^i, S_k^i),\,1\leq i\leq n,\, k=1,2,\dots\big)$.

\begin{lemma}
\label{fourthlemma}
Let $F$ and $F_m,\, m=1, 2, \dots$ 
be probability distribution
functions  with $F(0)=0$ and  $F_m(0)=0$ for all
$m=1,2,\dots$. 
If $F_m$ converges weakly to $F$ 
then $\mu(F_m; \rho_n, \ldots, \rho_1)$
also  converges weakly to $\mu(F; \rho_n, \ldots, \rho_1)$,
as $m\to\infty$,
\end{lemma}

%\noindent
%{\sl Remark:} 
%Note, that Lemma 1 and it's proof does not make sense for
%degenerate distributions  concentrated at $0$. This is why we
%have to impose the conditions on the distribution functions.

\begin{proof}[Sketch of proof]

The most natural (and rather standard) way to see this is by
use  of a space-time diagram.
This enables us to couple two on-off systems with the same
recovery rates but 
different input interval distributions, say $F$ and $F'$. We give 
a short outline of the argument:
Let $0 < I_1 < I_2 < I_3 , \ldots$ denote the points of a renewal
process 
with interval distribution $F$. (That is,  
$(I_{k+1} - I_k)_{k = 1, 2, \ldots}$ 
are i.i.d. random variables with distribution function $F$).
Now assign to each node $i$, independently of the other nodes and
of the above  
renewal process, a Poisson point process with intensity
$\rho_i$. These Poisson 
points are interpreted as \emph{potential} recovery points. This
means that 
if $t$ is such a point for node $i$, and node $i$  is in state
$0$ just before 
time $t$, it switches to state 1 at time $t$ (otherwise the point
is ignored). 
The $R_k^i$ and $S_k^i$ can be defined in a natural way in terms
of the above 
Poisson processes and the renewal process. If we now replace
$F$ by $F'$, we can compare the new situation with the old one
with the help of a suitable natural coupling: 
use the same realization of the above mentioned Poisson point
processes and 
take an obvious coupling of $F$ and $F'$. Details are left to the
reader.

\end{proof}

We will need the following notation: If $F$ is the input interval
distribution function at node $n$, then 
let, for $i \leq n$,
$F_{(\rho_n, \cdots, \rho_i)}$ denote the probability
distribution function 
of the intervals between successive signals received at node 
$i$, i.e. the  distribution of $S_1^i$. 
It is clear from the description of the system that for $i\le k\le n$ 
$$
\big(F_{(\rho_n,\dots,\rho_k)}\big)_{(\rho_{k-1},\ldots,\rho_i)} 
= 
F_{(\rho_n, \ldots, \rho_i)}.
$$
If $F$ and $G$ are two probability distribution functions, we write
$F \preceq G$ (or $G \succeq F$) if for any $x$ we have
$F(x) \geq G(x)$, i.e., if the
distribution $G$ stochastically dominates the distribution $F$.  
We have the following lemma:

\begin{lemma}
\label{fifthlemma}
For any $ \rho_k, \ldots, \rho_n, \rho_{n+1} > 0$ and any
probability  distribution function $F$,
$$
F_{(\rho_{n+1}, \rho_n, \ldots, \rho_k)} 
\succeq 
F_{(\rho_n, \ldots, \rho_k)}.
$$
\end{lemma}

\begin{proof}
Using Lemma 3 we have 
$F_{(\rho_{n+1}, \rho_n, \ldots, \rho_k)} =
F_{(\rho_n, \ldots, \rho_k, \rho_{n+1})} =$ \\
$\big(F_{(\rho_n, \ldots, \rho_k)}\big)_{(\rho_{n+1})}$, 
which obviously stochastically dominates 
$F_{(\rho_n, \ldots, \rho_k)}$.
\end{proof}

\noindent
{\bf Remarks:}
(i)  This lemma is not as obvious as it looks. For instance,
it is {\em not} true in general that if $F \succeq G$, then 
$F_{(\rho)} \succeq G_{(\rho)}$. The above argument
\emph{essentially relies} on Lemma 3.
\\
(ii)
Although, strictly speaking, Lemma 3 has not been proved for the 
case permanent  input (i.e. the case where the input interval
distribution function is $[0]$), its analog for 
that case follows easily from the fact that, as remarked earlier,
for such input signals
the sequence of signals received at node $n$ (and sent to node
$n-1$) corresponds to a
Poisson process 
with rate $\rho_n$, so that formally
\begin{equation}
\label{infeq}
[0]_{(\rho_n, \ldots, \rho_k)} = 
E^{\rho_n}_{(\rho_{n-1}, \ldots, \rho_k)},
\end{equation}
with $E^{\rho_n}$ the exponential distribution with mean
$1/\rho_n$.
In the sequel we shall use this notation for the exponential
distribution. 
Using (\ref{infeq}), and the (easy to check) fact that, if $E^{\rho}$ and
$E^{\rho'}$ are  exponential  distributions
with parameter $\rho$ and $\rho'$, respectively,  then 
$$
\big(E^{\rho}\big)_{(\rho')}=
E^{\rho}*E^{\rho'}=
\big(E^{\rho'})_{(\rho)},
$$ 
one can easily extend Lemma 3 to the case $F = [0]$.

The following lemma is a deterministic statement. First we give
some more  definitions and terminology.

A collection of non-negative numbers 
$s_{k-1}^i, r_k^i$, $1 \leq i \leq n$,\, $ k=1, 2, \ldots$ 
is called a (volume-$n$) \emph{signal/recovery} sequence if the
following hold:

\begin{enumerate}[(i)]
\item
For each $i$, $0= s_0^i < r_1^i < s_1^i < r_2^i < s_2^i <
\ldots$.  
\item
For each $i$ the set $\{s_{k-1}^i, r_k^i : k =1, 2, \ldots \}$ is
discrete. 
\item
For each $i < n$, and $k\ge 1$ 
$$
s^i_k=\min\{s^{i+1}_l: s^{i+1}_l>r^i_k\}.
$$
\end{enumerate}

The motivation for this definition is that the $r_k^i$-s and
$s_k^i$-s can be  interpreted as a realization of the recovery
and reception times  in an  on-off system.

We denote ${\cal S}^i:=\{s^i_k:k\ge0\}$. Property (iii) above is
clearly  equivalent to (iiia) \emph{and} (iiib) below, valid for
all $i<n$

{(iiia)}
${\cal S}^i\subset {\cal S}^{i+1}$

{(iiib)}
${\cal S}^{i+1}\setminus{\cal S}^{i}\subset
\cup_{k=1}^\infty(s^i_{k-1},r^i_{k}]$.

We now give a natural infinite version of this definition:
A collection of non-negative numbers 
$s_{k-1}^i, r_k^i$, $i = 1,2, \ldots$, $k=1,2, \ldots$ 
is called an infinite signal/recovery sequence if for
each $n$ the collection
$s_{k-1}^i, r_k^i$, $i=1,2\dots,n$, $k=1, 2, \ldots$ 
is a volume-$n$ signal/recovery sequence.
We say that the sequence has \emph{dense signals} if for every 
interval $V \subset \R_+$ there exist $i, k$ s.t. $s_k^i \in V$.
When $(s_{k-1}^i, r_k^i)$ is a, finite or infinite, signal/recovery
sequence, we define its corresponding {\em on-off sequence} as the following 
sequence of functions $\omega_i : R_+ \mapsto \{0,1\}, \, i \in \N.$
$$
\omega_i(t) := 
\left\{ 
\begin{array}{ll}
0      & \mbox{ \ if \ } t\in \cup_{k=1}^\infty
[s_{k-1}^i,r_{k}^i),
\\ [5pt]
1      & \mbox{ \ if \ } t\in \cup_{k=1}^\infty
[r_{k}^i,s_{k}^i)
\end{array}
\right. 
$$

\begin{lemma}
\label{sixthlemma}
Let $s_{k-1}^i,r_k^i$, $i=1, 2, \ldots, \, k=1, 2, \ldots$ 
be an infinite
signal/\-recovery sequence with dense signals. Let
$\omega_i(\cdot)$, $i=1,2\dots$  be the
corresponding on-off sequence. Then:
\begin{enumerate}[(a)]
\item
For each $i=1,2\dots$, the function $t\mapsto\omega_i(t)$ is 
c.a.d.l.a.g. 
\item
There are no $t$ and $k$ for which $\omega_l(t) = 1$, for all
$l\ge k$. 
\item
For each $t$ and $k$ with $\omega_l(t^-) = 1$ for all $l\ge k$,
we have $\omega_k(t) = 0$.
\item
For every $k$, \, $l > k$ and $t>0$ with 
$\omega_k(t^-) = 1$ and $\omega_l(t^-)= 0$, 
we have $\omega_k(t) = 1$.
\end{enumerate}
\end{lemma}

\begin{proof}

(a)
The c.a.d.l.a.g. property  
follows immediatley from the definition of the functions 
$\omega_i$.

(b) 
Suppose that for some $k$ and $t\ge0$ \, 
$\omega_l(t) = 1$ for all $l\ge k$.
Then, because of (a),
there is an $\vareps > 0$ such that  $\omega_k(s) = 1$ for all
$s\in (t, t+ \vareps)$.
Hence, by definition of $\omega_k$, 
${\cal S}^k \cap (t, t+\vareps) = \emptyset$.
However, because signals are dense, there is a $j> k$ with
${\cal S}^j \cap (t, t+ \vareps) \neq \emptyset$. 
Let $\tilde{j}$ be the smallest of such $j>k$.
So we have
$\omega_{\tilde j-1}(t)=1$, 
${\cal S}^{\tilde j-1}\cap(t,t+\vareps)=\emptyset$
and  
${\cal S}^{\tilde{j}} \cap (t, t+\vareps) \neq \emptyset$, 
which contradicts property (iiib) of a signal/recovery sequence. 

(c)
Suppose that for some $t>0$ and some $k$
$\omega_l(t^-) = 1$ for all $l\ge k$,  and
$\omega_k(t) = 1$.
% So for each $i \geq k$ there exists an $\vareps_i$ s.t. 
% $\omega_i(s) =1 1$ for  $s\in(t-\vareps_i, t)$. 
By (b) there is an $l>k$ with $\omega_l(t) = 0$.
Let $m$ be the smallest. So we have:
$\omega_{m-1}(t^-)=\omega_{m-1}(t)=1$
 and  
$\omega_m(t^-)=1,\,\omega_m(t)=0$.
This clearly implies that $t\in{\cal S}^{m}$ but at the same time
$t$ is in the interior of the set
$\cup_k[r^{m-1}_k,s^{m-1}_k)$. This contradicts property (iiib)
of signal/recovery systems. 

(d) 
Suppose $\omega_k(t^-) =1$ and $\omega_l(t^-) = 0$ for some $l>k$
and  $\omega_k(t) =0$. So $t \in {\cal S}^k$. But then (by property
(iii) of a signal/recovery system) $t \in {\cal S}^l$, which is
in  conflict with the above mentioned fact
that $\omega_l(t^-) = 0$.

\end{proof}

We continue our proof of Theorem 1. 
Let the $\rho_i$, $i=1, 2,\ldots$ be as in the statement of the
theorem, 
i.e., for all $t>0$,  $\sum_i e^{- \rho_i t} < \infty$.
Let, for $k \leq l$, $F^{(k,l)} = [0]_{(\rho_l, \ldots, \rho_k)}$
(see  earlier in this section). Using Lemma 5 we have
$$
F^{(k,l+1)} = [0]_{(\rho_{l+1}, \rho_l,\ldots, \rho_k)} 
\succeq
[0]_{(\rho_{l}, \ldots, \rho_k)} = F^{(k,l)}.
$$
Hence, keeping $k$ fixed,  the sequence of distributions
$F^{(k,l)}$, $l\ge k$,  converges weakly, as $l \to \infty$. The
following 
lemma shows that it converges to a probability distribution:

\begin{lemma}
\label{seventhlemma}
For each $k$, 
$F^{(k,l)}(t) \to 1$ as $t\to\infty$, uniformly in $l$.
\end{lemma}

\begin{proof}

As before, let $E^{\rho}$ 
denote the exponential distribution with mean $1/\rho$. 
For each $\rho > 0$ and $t > 0$ we have (using Lemma 5 again)
\begin{equation}
F^{(k,l)}(t)  
=  
[0]_{(\rho_l, \ldots, \rho_k)}(t) 
\preceq
[0]_{(\rho, \rho_l, \ldots, \rho_k)}(t) 
= 
E^{\rho}_{(\rho_l, \ldots,\rho_k)}(t).
\end{equation}
Note that this last expression is the probability that in a
finite on-off 
system with $l-k+1$ nodes with recovery rates $\rho_l, \ldots, \rho_k$,
and where the input signals are generated according to a Poisson process
with intensity $\rho$, the last node receives a signal before time $t$.
This probability is clearly larger than or equal to the probability 
that each of (a)-(c) below happens: \\
(a) No input signal is sent in the interval $(0, \sqrt t)$. \\
(b) Every node is in state 1 at time $\sqrt t$.\\
(c) An input signal is sent in the interval $(\sqrt t, t)$. \\
This probability is 
$$ 
e^{-\rho\sqrt{t}}
\prod_{j=k}^l(1-e^{-\rho_jt})
(1-e^{-\rho(t-\sqrt{t})})
\le
e^{-\rho \sqrt t} 
(1 - \sum_{j= 1}^\infty e^{ - \rho_j \sqrt t})
(1 - e^{- \rho (t - \sqrt t)}).
$$
For every $\rho$ this is a lower bound for $F^{(k,l)}(t)$. 
Now use (\ref{casefour}) and take
$\rho = t^{-2/3}$ to complete the proof of Lemma 7.

\end{proof}

We go on with the proof of Theorem 1.
We have seen that $F^{(k,l)}$ converges to a probability distribution
function as $l \to \infty$. Denote the limit by $F^{(k)}$, and let
$$ 
\mu_k := \mu(F^{(k)}; \rho_{k-1}, \ldots, \rho_1),
$$
where we use the notation introduced at 
the beginning of this section.
In this way we get a sequence $(\mu_k)$ of probability
measures on $\Sigma^{k-1}$,
where  $\Sigma$ is the set of all sequences 
$(s_{k-1},r_k)_{k=1}^{\infty}$
with $0 =s_0 < r_1 < s_1 < r_2 < \ldots$.
 From the definitions it is clear that for each $l$, the
projection of 
$\mu(F^{(k+1,l)}; \rho_k, \ldots, \rho_1)$ 
on $\Sigma^{k-1}$ equals
$\mu(F^{(k,l)}; \rho_{k-1}, \ldots, \rho_1)$.
By Lemma 4 it follows that the projection
of $\mu_{k+1}$ on $\Sigma^{k-1}$ is $\mu_k$. Hence, by standard
extension theorems, there is a measure $\nu$ on $\Sigma^{\N}$ whose
marginal on $\Sigma^k$ is $\mu_k$, $k=1, 2, \ldots$.
It is clear that for each $k$ a random element of $\Sigma^k$ is
$\mu_k$-a.s. a (volume $k$) signal/recovery sequence. Hence, a random
element of $\Sigma^{\N}$ is $\nu$-a.s. an infinite signal/recovery
sequence.
The theorem now follows from Lemma 6
if we can show that $\nu$-a.s. the system has dense signals.
By standard countability arguments this is equivalent to showing
that for every open interval $I \subset \R_+$,
\begin{equation}
\label{dense}
\nu \big\{\exists k: {\cal S}^k \cap  I \not=\emptyset \big\} = 1.
\end{equation}
Due to property (iiia) of signal/recovery systems, the l.h.s. of 
(\ref{dense}) equals 
$\lim_{k \to \infty}
\nu \big\{ {\cal S}^k \cap  I \not=\emptyset \big\}$
which, by the 
construction of $\nu$ above, equals
$$
\lim_{k \to \infty} \lim_{l \to \infty} 
\mu([0]; \rho_l, \ldots, \rho_k) 
\big\{ {\cal S}^k \cap  I \not=\emptyset \big\}.
$$

The required result now follows from the following Lemma:

\begin{lemma}
\label{eighthlemma}
For every open interval 
$I \subset {\R_+}$ 
and for every $\varepsilon > 0$
there exists a finite $K$ such that for all 
$k \geq K$ and $l \geq k$
$$
\mu([0]; \rho_l, \ldots, \rho_k) 
\big\{ {\cal S}^k \cap  I \not=\emptyset \big\} 
> 1- \varepsilon.
$$
\end{lemma}

\begin{proof}

We have, for any $\rho > 0$,
\begin{align}
\mu([0]; \rho_l, \ldots, \rho_k) 
\big\{  {\cal S}^k \cap  I \not=\emptyset \big \}
& 
\geq  
\mu([0]; \rho_l, \ldots, \rho_k, \rho)
\big \{  {\cal S}^{k-1} \cap  I \not=\emptyset  \big\} 
\notag
\\
&
= 
\mu([0]; \rho, \rho_l,  \ldots, \rho_k) 
\big\{ {\cal S}^k \cap  I \not=\emptyset  \big\} 
\notag 
\\
&
= 
\mu(E^{\rho}; \rho_l,  \ldots, \rho_k) 
\big\{ {\cal S}^k \cap  I \not=\emptyset   \big\},
\label{dlemma}
\end{align}
where the first two  expressions in the r.h.s. refer to a system with
leftmost and rightmost nodes $k-1$ and $l$, and $k$ and $l+1$, 
respectively.
The inequality is obvious from the definition, the first equality
follows from Lemma 3. 
Remind that $E^\rho$ denotes the exponential distribution function
with mean $1/\rho$. 
Note that the last expression in the r.h.s. of (\ref{dlemma})
is the probability that in a (size $l-k+1$)
on-off system
to which input signals are sent according to a Poisson process
with intensity 
$\rho$, and with recovery rates $\rho_l, \ldots, \rho_k$, the
last node receives 
a signal in the time interval $I$, and the computations below refer to that
system. We will choose 
$\rho$ appropriately, depending on $k$.
First of all, it follows from (\ref{casefour})
that there exists a sequence $(\tau_i)$ with the properties that
$\lim_{i \to \infty} \tau_i =0$, $\tau_i < |I|/2$ for all $i$, and
$\lim_{i\to\infty}\sum_{j \geq i} e^{-\rho_j \tau_i}= 0$. 
Now take $\rho = 1/ \sqrt{\tau_k}$.
Let $t$ and $t+s$ be the infimum and supremum of the interval $I$.
It is clear that the last expression in (\ref{dlemma}) is larger than or
equal to
the probability that each of the following events (a)-(c) occur:
\begin{enumerate}[(a)]
\item
No input signal is sent in $(t, t+\tau_k)$. 
\item
Each node in the system which had value $0$ at time $t$, has recovered
before time $t+\tau_k$.
\item
An input signal is sent in the interval $(t+s/2, t+s)$.
\end{enumerate}

\noindent
This probability is
$$
e^{-\sqrt{\tau_k}}
\prod_{j=k}^l\big(1-e^{-\rho_j\tau_k}\big)
\big(1-e^{-\frac{s}{2\sqrt{\tau_k}}}\big)
\ge
e^{-\sqrt{\tau_k}}
\big(1-\sum_{j\ge k}e^{-\rho_j\tau_k}\big)
\big(1-e^{-\frac{s}{2\sqrt{\tau_k}}}\big).
$$
The right hand side in the last inequality 
does not depend on $l$ and goes to $1$ as 
$k \to\infty$.
This completes the proof of Lemma 8 and of Theorem 1. 

\end{proof}

\subsubsection*{Acknowledgments}
The cooperation between the authors, which led to this paper, was
financially supported
by a Netherlands-Hungarian cooperation grant from the respective 
national research foundations NWO and OTKA.\\
JvdB had
interesting discussions with R. Meester and A. Luk\'acs.
The latter also drew our attention to the paper by Flajolet and
Sedgewick. \\
We also thank I. Benjamini and O. Schramm for allowing us to
mention their non-existence result for `frozen percolation' on
the square lattice, and A. J\'arai for allowing us to state
his 1-d non-existence example.

\bigskip
\noindent
------------------------------

\bigskip
\hbox{
\vbox{\hsize=5cm\noindent
{\sc
Jacob van den Berg
\\
CWI
\\
Kruislaan 413
\\
1098 SJ Amsterdam
\\
The Netherlands
\\
}}
\hskip2cm
\vbox{\hsize=5cm\noindent
{\sc 
B\'alint T\'oth
\\
Institute of Mathematics
\\
TU Budapest
\\
Egry J\'ozsef u. 1.
\\
H-1111 Budapest
\\
Hungary}}
}


\begin{thebibliography}{9}

\bibitem{aldous}
Aldous, D. (2000), 
The Percolation Process on a Tree where Infinite Clusters are Frozen.
{\sl Math. Proc. Cambridge Philos. Soc.} {\bf 128}, 465--477. 


\bibitem{benjaminischramm}
Benjamini, I., Schramm, O. (1999),
private communication (via D. Aldous).

\bibitem{flajoletsedgewick}
Flajolet, P., Sedgewick, R. (1995),
Mellin transforms and asymptotics: Finite differences and Rice's
integrals. 
{\sl Theoret. Comput. Sci.} {\bf 144}: 101--124.
 
\bibitem{jarai}
J\'arai, A. (1999),
private communication

\bibitem{lukacs}
Luk\'acs, A. (1999),
private communication

\end{thebibliography}
\end{document}